\documentclass[11pt,a4paper]{amsart}
\usepackage{amsmath,amssymb,amsfonts,mathtools}
\usepackage[margin=28mm]{geometry}
\usepackage[hidelinks]{hyperref}
\allowdisplaybreaks[2]
\newtheorem{theorem}{Theorem}[section]
\newtheorem{lemma}[theorem]{Lemma}

\newtheorem{corollary}[theorem]{Corollary}
\newtheorem{quest}[theorem]{Question}
\theoremstyle{remark}
\newtheorem{remark}[theorem]{Remark}
\numberwithin{equation}{section}
\newcommand{\M}{\mathcal M}
\newcommand{\N}{\mathcal N}
\newcommand{\A}{\mathcal A}
\newcommand{\R}{\mathcal R}

\newcommand{\bbN}{\mathbb N}
\newcommand{\bbT}{\mathbb T}
\newcommand{\bbM}{\mathbb M}
\newcommand{\E}{\mathbb E}
\newcommand{\Prob}{\mathbb P}
\newcommand{\Tr}{\operatorname{Tr}}
\newcommand{\tr}{\operatorname{tr}}
\newcommand{\supp}{\operatorname{supp}}

\newcommand{\tensor}{\mathbin{\bar\otimes}}
\newcommand{\norm}[1]{\left\|#1\right\|}

\newcommand{\emb}{\hookrightarrow}
\newcommand{\nemb}{\not\hookrightarrow}
\newcommand{\ind}{\mathbf 1}
\title[Non-equivalence of $L_p$-embeddings]{Non-equivalence of $L_p$-embeddings of symmetric operator ideals and their commutative cores}
\author{Chenhao Gao}
\author{Jinghao Huang}
\thanks{This work was supported by the National Natural Science Foundation of China (Nos. 12301160, 12471134 and 12671159).}
\begin{document}
\begin{abstract}
For every $1\leq p<2$, we construct a separable symmetric sequence space $E$ such that the corresponding operator ideal $C_E$ is isomorphic to a subspace of $L_p(\M,\sigma)$ for some finite von Neumann algebra $(\M,\sigma)$, whereas $E$ is not isomorphic to a subspace of $L_p(0,1)$. This answers a question (see Question~\ref{JSZquestion}) of Junge, Sukochev and Zanin in the negative.
\end{abstract}
\maketitle

\section{Introduction}

Let $H$ be a separable infinite-dimensional Hilbert space, and let $B(H)$ and $K(H)$ denote the algebras of all bounded and all compact linear operators on $H$, respectively. Given a symmetric Banach sequence space $E\subset c_0$, we denote by $C_E$ the symmetric operator ideal
\[
 C_E=\{A\in K(H):\mu(A)\in E\},\qquad
 \norm{A}_{C_E}=\norm{\mu(A)}_E,
\]
where $\mu(A)=(\mu_k(A))_{k\geq1}$ is the sequence of singular values of $A$ arranged in non-increasing order. In particular, $C_p=C_{\ell_p}$ is the Schatten--von Neumann ideal. The correspondence between symmetric sequence spaces and symmetric operator ideals goes back to Calkin \cite{Calkin}; the assertion that $C_E$ is a Banach ideal for an arbitrary symmetric Banach sequence space $E$ is due to Kalton and Sukochev \cite[Theorem~8.11]{KS}. Fixing an orthonormal basis in $H$, we identify $\ell_\infty$ with the algebra of diagonal operators in $B(H)$. The space $C_E\cap\ell_\infty$, equipped with the norm induced by $C_E$, is called the commutative core of $C_E$. Under the diagonal identification, the commutative core of $C_E$ coincides isometrically with $E$.

Throughout the paper, $L_p(\M,\tau)$ denotes the noncommutative $L_p$-space associated with a von Neumann algebra $\M$ equipped with a faithful normal semifinite trace $\tau$. For $1\leq p<\infty$, the norm of $x\in L_p(\M,\tau)$ is given by
\[
 \norm{x}_{L_p(\M,\tau)}=\tau(|x|^p)^{1/p}.
\]
We refer to \cite{PX} for the theory of noncommutative $L_p$-spaces. We say that $(\M,\tau)$ is a noncommutative probability space if $\tau(1)=1$. For Banach spaces $X$ and $Y$, we write $X\emb Y$ if there exists an isomorphic embedding of $X$ into $Y$, and $X\nemb Y$ otherwise. We write $X\simeq Y$ if $X$ and $Y$ are isomorphic as Banach spaces.

Let $1\leq p<2$. In \cite[Theorem~3]{JSZ}, Junge, Sukochev and Zanin established that if $E$ is a separable symmetric sequence space, $E\neq\ell_p$ (in the sense of sets) and $E\emb L_p(0,1)$, then $C_E\emb L_p(\R)$, where $\R$ is the hyperfinite $\mathrm{II}_1$ factor equipped with its normalized trace. It was also proved in \cite[Theorem~6]{JSZ} that if $E$ is a separable symmetric sequence space and $E\neq\ell_p$, then the following conditions are equivalent:
\begin{enumerate}
\item[(i)] $C_E\emb L_p(\M,\sigma)$ for some noncommutative probability space $(\M,\sigma)$;
\item[(ii)] $E\emb L_p(\N,\tau)$ for some noncommutative probability space $(\N,\tau)$.
\end{enumerate}
A continuous counterpart of the latter result, as well as an extension to general semifinite von Neumann algebras, was subsequently established in \cite{HJSZ}. On the other hand, the question whether condition (i) implies that $E\emb L_p(0,1)$ was left open in \cite[Question~8]{JSZ}, i.e.
\begin{quest}\label{JSZquestion}\cite[Question 8]{JSZ}
			Suppose that  (an ideal of compact operators) $C_E$ admits an isomorphic embedding into $\mathcal{L}_r(\M,\sigma)$, $1\le r<2$, for some finite von Neumann algebra $(\M,\sigma)$. Is it true that the commutative core $E$ of the ideal 
			$C_E$ admits an isomorphic embedding into $L_r(0,1)$? 
		\end{quest} 

The aim of the present paper is to answer Question~\ref{JSZquestion} in the negative. More precisely, we establish the following result.

\begin{theorem}\label{thm:main}
For every $1\leq p<2$, there exist a separable symmetric Banach sequence space $E$ and a noncommutative probability space $(\M,\sigma)$ such that
\[
 C_E\emb L_p(\M,\sigma),\qquad E\nemb L_p(0,1).
\]
\end{theorem}

By \cite[Theorem~6]{JSZ}, in order to prove Theorem~\ref{thm:main}, it suffices to construct a separable symmetric sequence space $E$ such that $E\emb L_p(\N,\tau)$ for some noncommutative probability space $(\N,\tau)$ and $E\nemb L_p(0,1)$. There are numerous examples of subspaces of $L_p(\N,\tau)$ with an unconditional basis which do not embed into $L_p(0,1)$ (see, in particular, \cite[Corollary~1.4]{HRS}). However, the existence of such a subspace with a symmetric basis does not follow directly from \cite[Corollary~1.4]{HRS}.

The principal step in the proof of Theorem~\ref{thm:main} is the following result, which is a vector-valued version of the symmetrization theorem of Maurey and Schechtman \cite[Main Theorem]{MS}. For earlier results concerning the embedding of spaces with unconditional bases into spaces with symmetric bases, see \cite{Davis,Lindenstrauss,Szankowski}.

\begin{theorem}\label{thm:symm}
Let $1\leq p<2$, and let $X$ be an infinite-dimensional Banach space with an unconditional Schauder basis and cotype $2$. There exist a separable symmetric Banach sequence space $E$, a partition $(B_i)_{i\geq1}$ of $\bbN$ into successive finite intervals and a sequence of positive scalars $(b_i)_{i\geq1}$ such that
\begin{enumerate}
\item[(i)] $E\emb L_p(0,1;X)$;
\item[(ii)] the closed linear span of the sequence
\[
 u_i=b_i\sum_{k\in B_i}e_k,\qquad i\geq1,
\]
is isomorphic to $X$ and is the range of a contractive projection on $E$.
\end{enumerate}
Here $(e_k)_{k\geq1}$ denotes the unit vector basis of $E$.
\end{theorem}

In contrast to the classical symmetrization theorem \cite{MS}, the Banach space $X$ in Theorem~\ref{thm:symm} need not embed into a commutative $L_p$-space. This feature is essential for the application to Theorem~\ref{thm:main}.

\section{Symmetrization of unconditional bases}

We work over the complex scalar field. For a probability space $(\Omega,\Prob)$ and a Banach space $X$, $L_p(\Omega;X)$ denotes the Bochner space of strongly measurable $X$-valued functions with finite norm
\[
 \norm{f}_{L_p(\Omega;X)}=\big(\E\norm{f(\omega)}_X^p\big)^{1/p}.
\]
A Steinhaus random variable is a complex-valued random variable distributed according to normalized Haar measure on $\bbT=\{z\in\mathbb C:|z|=1\}$.

Recall that a Banach space $X$ is of cotype $2$ if there exists $K\geq1$ such that, for every $m\geq1$ and $v_1,\ldots,v_m\in X$,
\begin{equation}\label{eq:cotype}
 \left(\sum_{j=1}^m\norm{v_j}_X^2\right)^{1/2}
 \leq K\left(\E_\varepsilon\norm{\sum_{j=1}^m\varepsilon_j v_j}_X^2\right)^{1/2},
\end{equation}
where $(\varepsilon_j)_{j=1}^m$ is a sequence of independent Rademacher random variables. It is well known that $L_p(\M,\tau)$ is of cotype $2$ for $1\leq p\leq2$ (see \cite{Fack,TJ} and \cite[Theorem~13]{JSZ}). In particular, every subspace of $L_p(\M,\tau)$ is of cotype $2$.

If $(x_i)_{i\geq1}$ is a normalized $1$-unconditional basis of $X$, then
\begin{equation}\label{eq:unconditional}
 \norm{\sum_i \lambda_i a_i x_i}_X\leq\norm{\sum_i a_i x_i}_X,
 \qquad |\lambda_i|\leq1,\quad (a_i)_{i\geq1}\in c_{00}.
\end{equation}
Equality holds in \eqref{eq:unconditional} if $|\lambda_i|=1$ for every $i$. In particular, if $v_1,\ldots,v_m$ are finitely supported vectors with pairwise disjoint supports with respect to $(x_i)_{i\geq1}$, then \eqref{eq:cotype} yields
\begin{equation}\label{eq:disjoint}
 \left(\sum_{j=1}^m\norm{v_j}_X^2\right)^{1/2}
 \leq K\norm{\sum_{j=1}^m v_j}_X.
\end{equation}

\begin{lemma}\label{lem:selectors}
Let $1\leq p\leq2$. Let $X$ have a normalized $1$-unconditional basis $(x_i)_{i\geq1}$ satisfying \eqref{eq:disjoint}. Fix an integer $m\geq1$ and independent $\{0,1\}$-valued random variables $(\xi_i)_{i\geq1}$ such that
\[
 \Prob(\xi_i=1)\leq\frac1m,\qquad i\geq1.
\]
For every $(a_i)_{i\geq1}\in c_{00}$, we have
\begin{equation}\label{eq:selectors}
 \left(\E\norm{\sum_i a_i\xi_i x_i}_X^p\right)^{1/p}
 \leq\frac{K}{\sqrt m}\norm{\sum_i a_i x_i}_X.
\end{equation}
\end{lemma}
\begin{proof}
Fix $(a_i)_{i\geq1}\in c_{00}$, and put $q_i=\Prob(\xi_i=1)$. Let $(t_i)_{i\geq1}$ be independent random variables uniformly distributed on $(0,1)$. The joint distribution of $(\xi_i)_{i\geq1}$ coincides with that of $(\ind_{\{t_i\leq q_i\}})_{i\geq1}$. Since $q_i\leq1/m$, \eqref{eq:unconditional} implies
\[
 \E\norm{\sum_i a_i\xi_i x_i}_X^p
 \leq \E\norm{\sum_i a_i\ind_{\{t_i\leq1/m\}}x_i}_X^p.
\]
Thus, it suffices to prove \eqref{eq:selectors} for independent selectors of parameter $1/m$.

Let $(\lambda_i)_{i\geq1}$ be independent random variables uniformly distributed on $\{1,\ldots,m\}$, and set
\[
 v_j=\sum_i a_i\ind_{\{\lambda_i=j\}}x_i,\qquad 1\leq j\leq m.
\]
For every realization of $(\lambda_i)_{i\geq1}$, the vectors $v_j$, $1\leq j\leq m$, have pairwise disjoint supports and $\sum_{j=1}^m v_j=\sum_i a_ix_i$. Applying \eqref{eq:disjoint}, taking expectations and using the fact that $v_1,\ldots,v_m$ are identically distributed, we obtain
\[
 m\E\norm{v_1}_X^2
 =\E\sum_{j=1}^m\norm{v_j}_X^2
 \leq K^2\norm{\sum_i a_ix_i}_X^2.
\]
Since $p\leq2$ and $(\ind_{\{\lambda_i=1\}})_{i\geq1}$ consists of independent selectors of parameter $1/m$, we have
\[
 \left(\E\norm{v_1}_X^p\right)^{1/p}
 \leq\left(\E\norm{v_1}_X^2\right)^{1/2}
 \leq\frac{K}{\sqrt m}\norm{\sum_i a_ix_i}_X.
\]
This completes the proof of the lemma.
\end{proof}

\begin{proof}[Proof of Theorem~\ref{thm:symm}]
Normalize the given unconditional basis $(x_i)_{i\geq1}$ of $X$. Replacing the norm of $X$ by the equivalent norm given on finitely supported vectors by
\[
 \left\|\sum_i a_ix_i\right\|_{\mathrm{u}}
 =\sup_{|\lambda_i|\leq1}\norm{\sum_i\lambda_i a_ix_i}_X,
\]
we may assume that $(x_i)_{i\geq1}$ is normalized and $1$-unconditional. This renorming preserves cotype and the isomorphism class of $L_p(0,1;X)$. Fix $K\geq1$ such that \eqref{eq:disjoint} holds. Set
\[
 \theta=\frac12\left(\frac12+\frac1p\right),\qquad
 \delta=\theta-\frac12=\frac1p-\theta=\frac{2-p}{4p}>0.
\]
Choose a positive integer $L$ and set $m_i=2^{Li^2}$, $i\geq1$, so that
\begin{equation}\label{eq:lacunary}
 \sum_{i=1}^\infty m_i^{-1}<1,\qquad
 2\sum_{1\leq i<j<\infty}\left(\frac{m_i}{m_j}\right)^\delta<\frac1{2e}.
\end{equation}
Indeed, writing $r=2^{-L\delta}$ and $j=i+h$ with $h\geq1$, we have
\[
 \sum_{i<j}\left(\frac{m_i}{m_j}\right)^\delta
 =\sum_{i,h\geq1}r^{2ih+h^2}
 \leq\sum_{i,h\geq1}r^{i+h}=\frac{r^2}{(1-r)^2}.
\]
The last expression tends to zero as $L\to\infty$. The first inequality in \eqref{eq:lacunary} also holds for all sufficiently large $L$.

On the Lebesgue probability space $(0,1)$, choose an independent sequence $(J_k)_{k\geq1}$ of identically distributed $\{0\}\cup\bbN$-valued random variables such that
\[
 \Prob(J_k=j)=\frac1{m_j},\quad j\geq1,\qquad
 \Prob(J_k=0)=1-\sum_{j=1}^\infty\frac1{m_j}.
\]
Choose a sequence $(\zeta_k)_{k\geq1}$ of independent Steinhaus random variables which is independent of $(J_k)_{k\geq1}$. Such sequences exist on $(0,1)$, since the countable product of the distribution space of $J_1$ and $\bbT$, equipped with the corresponding product probability measure, is a standard atomless probability space.

Define
\begin{equation}\label{eq:Z}
 Z_k=\sum_{j=1}^\infty m_j^\theta\zeta_k\ind_{\{J_k=j\}}x_j,\qquad k\geq1.
\end{equation}
At every point of the probability space, at most one summand in \eqref{eq:Z} is nonzero. The functions $Z_k$ are strongly measurable. Since $\theta p<1$ and $\norm{x_j}_X=1$, we have
\begin{equation}\label{eq:alpha}
 0<\alpha:=\norm{Z_1}_{L_p(0,1;X)}
 =\left(\sum_{j=1}^\infty m_j^{\theta p-1}\right)^{1/p}<\infty.
\end{equation}
In particular, every $Z_k$ belongs to $L_p(0,1;X)$.

For $a=(a_k)_{k\geq1}\in c_{00}$, set
\begin{equation}\label{eq:E}
 \norm{a}_E=\alpha^{-1}\norm{\sum_k a_k Z_k}_{L_p(0,1;X)}.
\end{equation}
We verify that \eqref{eq:E} defines a symmetric norm. The $X$-valued random variables $(Z_k)_{k\geq1}$ are independent and identically distributed, and $Z_k$ and $\omega Z_k$ have the same distribution for every $\omega\in\bbT$. Consequently, for every permutation $\pi$ of $\bbN$ and every sequence $(\omega_k)_{k\geq1}\subset\bbT$,
\begin{equation}\label{eq:symmetry}
 \norm{\sum_k a_ke_k}_E=\norm{\sum_k\omega_k a_ke_{\pi(k)}}_E,
 \qquad (a_k)_{k\geq1}\in c_{00}.
\end{equation}
Here $e_k$ is the $k$th unit vector of $c_{00}$. By averaging over the phases in \eqref{eq:symmetry}, we obtain
\begin{equation}\label{eq:coordinate}
 \norm{\sum_k d_k a_ke_k}_E\leq\norm{\sum_k a_ke_k}_E,
 \qquad |d_k|\leq1,\quad (a_k)_{k\geq1}\in c_{00}.
\end{equation}
Taking $d_k=\ind_{\{k=j\}}$ and using $\norm{e_j}_E=1$, we deduce that $\norm{a}_\infty\leq\norm{a}_E$. This proves positive definiteness of \eqref{eq:E}; homogeneity and the triangle inequality follow from the corresponding properties of the Bochner norm.

Let $E$ be the completion of $c_{00}$ in the norm \eqref{eq:E}. The coordinate projections $P_n(a)=\sum_{k=1}^n a_ke_k$, $n\geq1$, extend to contractions on $E$ by \eqref{eq:coordinate}, and $P_n a\to a$ in $E$ for every $a\in E$ by density of $c_{00}$. In view of the inequality $\norm{a}_\infty\leq\norm{a}_E$ on $c_{00}$, the completion embeds continuously into $c_0$: the continuous extension of the inclusion $c_{00}\subset c_0$ is injective because the projections $P_n$ converge strongly to the identity on $E$.

The completion is an ideal in $\ell_\infty$ by \eqref{eq:coordinate}. To verify rearrangement invariance, let $a\in E$ and $b\in c_0$ satisfy $\mu(a)=\mu(b)$, where $\mu(a)$ denotes the non-increasing rearrangement of $(|a_k|)_{k\geq1}$. Enumerate the nonzero coordinates of $a$ and $b$ in matching order of their absolute values. Denote by $a^{(n)}$ and $b^{(n)}$ the vectors obtained by retaining the first $n$ of the enumerated coordinates of $a$ and $b$, respectively (retaining all of them if there are fewer than $n$). By \eqref{eq:coordinate} and the density of $c_{00}$, $a^{(n)}\to a$ in $E$. By \eqref{eq:symmetry},
\[
 \norm{b^{(n)}-b^{(m)}}_E=\norm{a^{(n)}-a^{(m)}}_E,
 \qquad n,m\geq1.
\]
It follows that $(b^{(n)})_{n\geq1}$ converges in $E$ to $b$ and that $\norm{b}_E=\norm{a}_E$. Thus $E$ is a separable symmetric Banach sequence space. By \eqref{eq:E}, the linear mapping $S\colon E\to L_p(0,1;X)$ defined on the basis vectors by
\[
 Se_k=\alpha^{-1}Z_k,\qquad k\geq1,
\]
is an isometric embedding. This proves assertion (i).

Partition $\bbN$ into successive finite intervals $B_i$ of cardinality $m_i$, $i\geq1$, and set
\begin{equation}\label{eq:Yi}
 u_i=\alpha m_i^{-\theta}\sum_{k\in B_i}e_k,\qquad
 Y_i=Su_i=m_i^{-\theta}\sum_{k\in B_i}Z_k,\qquad i\geq1.
\end{equation}
For $i,j\geq1$, define the scalar random variables
\[
 \eta_{ij}=\sum_{k\in B_i}\zeta_k\ind_{\{J_k=j\}},\qquad
 N_{ij}=\sum_{k\in B_i}\ind_{\{J_k=j\}}.
\]
The distribution of $N_{ij}$ is binomial with parameters $m_i$ and $1/m_j$. By \eqref{eq:Z} and \eqref{eq:Yi},
\begin{equation}\label{eq:decomp}
 Y_i=\eta_{ii}x_i+R_i,\qquad
 R_i=\sum_{\substack{j\geq1\\j\neq i}}
       \left(\frac{m_j}{m_i}\right)^\theta\eta_{ij}x_j.
\end{equation}
For each realization of $(J_k)_{k\in B_i}$, the sum defining $R_i$ has at most $m_i$ nonzero terms. The series in \eqref{eq:decomp} also converges in $L_p(0,1;X)$, as shown by the estimates below.

\medskip
\noindent\textit{Estimates of the sequence $(\eta_{ii}x_i)_{i\geq1}$.}
Write $\eta_i=\eta_{ii}$ and $N_i=N_{ii}$ for $i\geq1$. Since $N_i$ has binomial distribution with parameters $m_i$ and $1/m_i$, we have
\begin{equation}\label{eq:onehit}
 \E|\eta_i|\geq\Prob(N_i=1)
 =\left(1-\frac1{m_i}\right)^{m_i-1}\geq e^{-1}.
\end{equation}
The last inequality follows from $(m_i-1)\log(1+1/(m_i-1))\leq1$. By \eqref{eq:unconditional}, Jensen's inequality and \eqref{eq:onehit}, for every $(a_i)_{i\geq1}\in c_{00}$, we obtain
\begin{align}
 \norm{\sum_i a_i\eta_i x_i}_{L_p(0,1;X)}
 &\geq\E\norm{\sum_i a_i|\eta_i| x_i}_X\notag\\
 &\geq\norm{\sum_i a_i\E|\eta_i| x_i}_X
 \geq e^{-1}\norm{\sum_i a_ix_i}_X.\label{eq:lower}
\end{align}

For every integer $\ell\geq1$, using the convention that $\binom{m_i}{\ell}=0$ for $\ell>m_i$, we have
\begin{equation}\label{eq:binomial}
 \Prob(N_i\geq\ell)\leq\E\binom{N_i}{\ell}
 =\binom{m_i}{\ell}m_i^{-\ell}\leq\frac1{\ell!}.
\end{equation}
For fixed $\ell\geq1$, the random variables $(\ind_{\{N_i\geq\ell\}})_{i\geq1}$ are independent, since the sets $B_i$ are pairwise disjoint. By \eqref{eq:binomial} and Lemma~\ref{lem:selectors} applied with $m=\ell!$,
\begin{equation}\label{eq:level}
 \norm{\sum_i a_i\ind_{\{N_i\geq\ell\}}x_i}_{L_p(0,1;X)}
 \leq\frac{K}{\sqrt{\ell!}}\norm{\sum_i a_ix_i}_X,
 \qquad (a_i)_{i\geq1}\in c_{00}.
\end{equation}
Since $|\eta_i|\leq N_i$ and $N_i=\sum_{\ell=1}^\infty\ind_{\{N_i\geq\ell\}}$, combining \eqref{eq:unconditional} and \eqref{eq:level}, we obtain
\begin{align}
 \norm{\sum_i a_i\eta_i x_i}_{L_p(0,1;X)}
 &\leq\norm{\sum_i a_iN_i x_i}_{L_p(0,1;X)}\notag\\
 &\leq\sum_{\ell=1}^\infty
 \norm{\sum_i a_i\ind_{\{N_i\geq\ell\}}x_i}_{L_p(0,1;X)}\notag\\
 &\leq K\left(\sum_{\ell=1}^\infty\frac1{\sqrt{\ell!}}\right)
 \norm{\sum_i a_i x_i}_X.\label{eq:upper}
\end{align}
For convenience, set $C_0=K\sum_{\ell=1}^\infty(\ell!)^{-1/2}<\infty$.

\medskip
\noindent\textit{Estimates of the sequence $(R_i)_{i\geq1}$.}
We first estimate $\eta_{ij}$ for $j<i$. Since the random variables $\zeta_k$ are independent, of mean zero and independent of $(J_k)_{k\geq1}$, we have
\[
 \E|\eta_{ij}|^2
 =\sum_{k\in B_i}\Prob(J_k=j)=\frac{m_i}{m_j}.
\]
As $p<2$, it follows that
\begin{equation}\label{eq:smallj}
 \left(\frac{m_j}{m_i}\right)^\theta\norm{\eta_{ij}}_{L_p(0,1)}
 \leq\left(\frac{m_j}{m_i}\right)^{\theta-1/2}
 =\left(\frac{m_j}{m_i}\right)^\delta,\qquad j<i.
\end{equation}
For $j>i$, condition on $(J_k)_{k\in B_i}$. The conditional second moment of $\eta_{ij}$ equals $N_{ij}$. Therefore,
\[
 \E_\zeta|\eta_{ij}|^p\leq N_{ij}^{p/2}\leq N_{ij}.
\]
Here $\E_\zeta$ denotes expectation with respect to the Steinhaus variables, with $(J_k)_{k\in B_i}$ fixed; the last inequality uses $N_{ij}\in\{0,1,\ldots,m_i\}$ and $p\leq2$. Taking expectations with respect to $(J_k)_{k\in B_i}$, we obtain
\[
 \E|\eta_{ij}|^p\leq\E N_{ij}=\frac{m_i}{m_j}.
\]
Consequently,
\begin{equation}\label{eq:largej}
 \left(\frac{m_j}{m_i}\right)^\theta\norm{\eta_{ij}}_{L_p(0,1)}
 \leq\left(\frac{m_i}{m_j}\right)^{1/p-\theta}
 =\left(\frac{m_i}{m_j}\right)^\delta,\qquad j>i.
\end{equation}
Combining \eqref{eq:decomp}, \eqref{eq:smallj} and \eqref{eq:largej}, and using $\norm{x_j}_X=1$ for every $j\geq1$, we obtain
\[
 \norm{R_i}_{L_p(0,1;X)}
 \leq\sum_{j<i}\left(\frac{m_j}{m_i}\right)^\delta
       +\sum_{j>i}\left(\frac{m_i}{m_j}\right)^\delta.
\]
In particular, the series in \eqref{eq:decomp} is absolutely convergent in the Bochner norm. Summing over $i\geq1$ and applying \eqref{eq:lacunary}, we arrive at
\begin{equation}\label{eq:errorsum}
 \sum_{i=1}^\infty\norm{R_i}_{L_p(0,1;X)}
 \leq2\sum_{i<j}\left(\frac{m_i}{m_j}\right)^\delta<\frac1{2e}.
\end{equation}
Since $(x_i)_{i\geq1}$ is a normalized $1$-unconditional basis, we have $\sup_i|a_i|\leq\norm{\sum_i a_ix_i}_X$ for every $(a_i)_{i\geq1}\in c_{00}$. Hence, by \eqref{eq:errorsum},
\begin{equation}\label{eq:error}
 \norm{\sum_i a_iR_i}_{L_p(0,1;X)}
 \leq\sup_i|a_i|\sum_i\norm{R_i}_{L_p(0,1;X)}
 \leq\frac1{2e}\norm{\sum_i a_ix_i}_X.
\end{equation}

Combining \eqref{eq:decomp}, \eqref{eq:lower}, \eqref{eq:upper} and \eqref{eq:error}, we infer that
\begin{equation}\label{eq:equivalence}
 \frac1{2e}\norm{\sum_i a_ix_i}_X
 \leq\norm{\sum_i a_iY_i}_{L_p(0,1;X)}
 \leq\left(C_0+\frac1{2e}\right)\norm{\sum_i a_ix_i}_X,
 \qquad (a_i)_{i\geq1}\in c_{00}.
\end{equation}
Since $Su_i=Y_i$ and $S$ is an isometric embedding, \eqref{eq:equivalence} shows that the linear mapping $U\colon X\to E$ given by $Ux_i=u_i$, $i\geq1$, is an isomorphic embedding.

Finally, define a linear mapping $P\colon c_{00}\to c_{00}$ by
\[
 Pa=\sum_{i=1}^\infty\left(\frac1{m_i}\sum_{k\in B_i}a_k\right)
                      \sum_{k\in B_i}e_k,\qquad a=(a_k)_{k\geq1}\in c_{00}.
\]
For $a\in c_{00}$, only finitely many intervals $B_i$ meet $\supp(a)$. The vector $Pa$ is the average of the vectors obtained from $a$ by permuting coordinates separately in each of these intervals. By \eqref{eq:symmetry} and the triangle inequality, $\norm{Pa}_E\leq\norm{a}_E$. Thus $P$ extends to a contraction on $E$. Since $P^2=P$ and $P(c_{00})$ coincides with the linear span of $(u_i)_{i\geq1}$, the range of $P$ is precisely $U(X)$. This proves assertion (ii), and completes the proof of the theorem.
\end{proof}

\begin{corollary}\label{cor:transfer}
Let $1\leq p<2$, and let $(\A,\tau)$ be a noncommutative probability space. If an infinite-dimensional Banach space $X$ with an unconditional Schauder basis is isomorphic to a subspace of $L_p(\A,\tau)$, then there exists a separable symmetric Banach sequence space $E$ such that $X$ is isomorphic to a complemented subspace of $E$ and
\[
 E\emb L_p\big(L_\infty(0,1)\tensor\A,dm\otimes\tau\big),
\]
where $dm$ is Lebesgue measure on $(0,1)$.
\end{corollary}
\begin{proof}
Let $T\colon X\to L_p(\A,\tau)$ be an isomorphic embedding. Since $L_p(\A,\tau)$ is of cotype $2$, Theorem~\ref{thm:symm} applies to $X$. Moreover, the mapping $f\mapsto T\circ f$ is an isomorphic embedding from $L_p(0,1;X)$ into $L_p(0,1;L_p(\A,\tau))$. For pairwise disjoint measurable sets $D_1,\ldots,D_m\subset(0,1)$ and $a_1,\ldots,a_m\in L_p(\A,\tau)$, we have
\[
 \left\|\sum_{j=1}^m \ind_{D_j}\otimes a_j\right\|_{L_p(L_\infty(0,1)\tensor\A)}^p
 =\sum_{j=1}^m |D_j|\tau(|a_j|^p)
 =\left\|\sum_{j=1}^m \ind_{D_j}a_j\right\|_{L_p(0,1;L_p(\A,\tau))}^p,
\]
where $|D_j|$ is the Lebesgue measure of $D_j$. The density of simple functions in the Bochner space gives an isometric embedding of $L_p(0,1;L_p(\A,\tau))$ into $L_p(L_\infty(0,1)\tensor\A,dm\otimes\tau)$. The product trace is a faithful normal tracial state. The assertion follows.
\end{proof}

\section{Proof of Theorem~\ref{thm:main}}

\begin{proof}[Proof of Theorem~\ref{thm:main}]
Fix $1\leq p<2$. For $n\geq1$, let $\bbM_n$ denote the algebra of all $n\times n$ complex matrices equipped with the normalized trace $\tr_n=n^{-1}\Tr$, where $\Tr$ is the usual matrix trace. Let $(e_{ij}^{(n)})_{1\leq i,j\leq n}$ be the system of matrix units in $\bbM_n$. Set
\[
 \A_n=\bbM_n\tensor L_\infty(\bbT^{n^2}),\qquad
 \tau_n=\tr_n\otimes m_n^{\mathrm H},
\]
where $m_n^{\mathrm H}$ denotes the normalized Haar measure on $\bbT^{n^2}$. Write $z_{ij}$, $1\leq i,j\leq n$, for the coordinate functions on $\bbT^{n^2}$, so that $(z_{ij})_{1\leq i,j\leq n}$ is a family of independent Steinhaus random variables. Define
\[
 x_{ij}^{(n)}=n^{1/p}e_{ij}^{(n)}\otimes z_{ij},\qquad
 F_n=\operatorname{span}\{x_{ij}^{(n)}:1\leq i,j\leq n\}
 \subset L_p(\A_n,\tau_n).
\]
The normalization factor $n^{1/p}$ ensures that $\norm{x_{ij}^{(n)}}_{L_p(\A_n,\tau_n)}=1$. The invariance of the product Haar measure under coordinatewise multiplication by unimodular scalars implies that $(x_{ij}^{(n)})_{1\leq i,j\leq n}$ is a $1$-unconditional basis of $F_n$.

For every $1\leq i,j\leq n$ and every $a_1,\ldots,a_n\in\mathbb C$, we have
\begin{align}
 \norm{\sum_{k=1}^n a_k x_{ik}^{(n)}}_{L_p(\A_n,\tau_n)}
 &=\left(\sum_{k=1}^n|a_k|^2\right)^{1/2},\label{eq:row}\\
 \norm{\sum_{k=1}^n a_k x_{kj}^{(n)}}_{L_p(\A_n,\tau_n)}
 &=\left(\sum_{k=1}^n|a_k|^2\right)^{1/2}.\label{eq:col}
\end{align}
Indeed, for each $z\in\bbT^{n^2}$, the matrices $\sum_{k=1}^n a_kz_{ik}e_{ik}^{(n)}$ and $\sum_{k=1}^n a_kz_{kj}e_{kj}^{(n)}$ have at most one nonzero singular value, equal to $(\sum_{k=1}^n|a_k|^2)^{1/2}$. The normalization by $n^{1/p}$ cancels the factor $n^{-1/p}$ coming from $\tr_n$.

If $1\leq i_1<\cdots<i_m\leq n$ and $1\leq j_1<\cdots<j_m\leq n$, then the nonzero singular values, counted with multiplicity, of the matrix $\sum_{k=1}^m a_kz_{i_kj_k}e_{i_kj_k}^{(n)}$ are the nonzero numbers among $|a_1|,\ldots,|a_m|$. Hence
\begin{equation}\label{eq:diag}
 \norm{\sum_{k=1}^m a_k x_{i_kj_k}^{(n)}}_{L_p(\A_n,\tau_n)}
 =\left(\sum_{k=1}^m |a_k|^p\right)^{1/p},\qquad a_1,\ldots,a_m\in\mathbb C.
\end{equation}
By \eqref{eq:row}--\eqref{eq:diag}, the space $F_n$ is a $1$-$GC_p^n$-space in the sense of \cite[Section~1]{HRS}. In particular, by \cite[Corollary~1.4]{HRS}, there is no constant $D<\infty$ such that for every $n\geq1$ there exists an isomorphic embedding $V_n\colon F_n\to L_p(0,1)$ satisfying $\norm{V_n}\norm{V_n^{-1}}\leq D$. Here $V_n^{-1}$ is defined on $V_n(F_n)$.

Consider the Banach space
\begin{equation}\label{eq:X}
 X=\left(\bigoplus_{n=1}^\infty F_n\right)_{\ell_p},\qquad
 \norm{(f_n)_{n\geq1}}_X=\left(\sum_{n=1}^\infty\norm{f_n}_{L_p(\A_n,\tau_n)}^p\right)^{1/p}.
\end{equation}
Ordering the bases of the $F_n$ successively gives a normalized $1$-unconditional Schauder basis of $X$. Since every $F_n$ embeds isometrically into $X$, the result of Haagerup, Rosenthal and Sukochev just recalled implies
\begin{equation}\label{eq:Xnon}
 X\nemb L_p(0,1).
\end{equation}
Indeed, if $V\colon X\to L_p(0,1)$ were an isomorphic embedding, its restrictions to the copies of $F_n$ in $X$ would have distortions bounded by $\norm{V}\norm{V^{-1}}$, independently of $n$, contradicting \cite[Corollary~1.4]{HRS}.

Let
\[
 \A=\bigoplus_{n=1}^\infty\A_n,\qquad
 \tau((a_n)_{n\geq1})=\sum_{n=1}^\infty2^{-n}\tau_n(a_n),
\]
where the von Neumann algebra direct sum is taken in the $\ell_\infty$ sense. The functional $\tau$ is a faithful normal tracial state on $\A$. For every finitely supported sequence $(f_n)_{n\geq1}$ with $f_n\in F_n$, we have
\[
 \norm{(2^{n/p}f_n)_{n\geq1}}_{L_p(\A,\tau)}^p
 =\sum_{n=1}^\infty2^{-n}\tau_n(|2^{n/p}f_n|^p)
 =\sum_{n=1}^\infty\norm{f_n}_{L_p(\A_n,\tau_n)}^p.
\]
Therefore the mapping $(f_n)_{n\geq1}\mapsto(2^{n/p}f_n)_{n\geq1}$ extends to an isometric embedding of $X$ into $L_p(\A,\tau)$.

By Corollary~\ref{cor:transfer}, there exists a separable symmetric Banach sequence space $E$ containing an isomorphic copy of $X$ as a complemented subspace and admitting an isomorphic embedding into $L_p(\N,\rho)$, where
\[
 \N=L_\infty(0,1)\tensor\A,\qquad \rho=dm\otimes\tau.
\]
Since $X\emb E$, \eqref{eq:Xnon} yields
\begin{equation}\label{eq:Enon}
 E\nemb L_p(0,1).
\end{equation}
In particular, $E\not\simeq\ell_p$, since $\ell_p$ embeds isometrically into $L_p(0,1)$ by means of normalized characteristic functions of pairwise disjoint measurable sets. Thus $E\neq\ell_p$ in the sense of sets as well: if the two spaces coincided as sets, their norms would be equivalent by the closed graph theorem, since both norms dominate the supremum norm.

Let $C_E$ be the symmetric operator ideal corresponding to $E$. The finite-rank operators are dense in $C_E$. Indeed, if $A\in C_E$ and $A_n$ is the truncation of the singular value decomposition of $A$ to its first $n$ singular values, then
\[
 \norm{A-A_n}_{C_E}
 =\norm{(\mu_k(A)\ind_{\{k>n\}})_{k\geq1}}_E\longrightarrow0,
\]
since $c_{00}$ is dense in $E$. For each fixed $n$, the set of operators of rank at most $n$ is separable in the $C_E$-norm: for an operator $B$ of rank at most $2n$, one has
\[
 \norm{B}_{C_E}\leq\norm{B}_{B(H)}\norm{\sum_{k=1}^{2n}e_k}_E,
\]
and the operators of rank at most $n$ form a separable metric space in the operator norm because $H$ is separable. Consequently, $C_E$ is a separable Banach ideal. Moreover, $C_E\neq C_p$, since their commutative cores are distinct.

We have verified all the hypotheses of \cite[Theorem~6]{JSZ}. Applying the implication (ii)$\Rightarrow$(i) of that theorem, we obtain a noncommutative probability space $(\M,\sigma)$ such that
\[
 C_E\emb L_p(\M,\sigma).
\]
Together with \eqref{eq:Enon}, this establishes Theorem~\ref{thm:main}.
\end{proof}

\begin{remark}
The von Neumann algebra $\N$ constructed above is hyperfinite and has separable predual. More explicitly, $\N$ is the direct sum of the algebras $\bbM_n\tensor L_\infty((0,1)\times\bbT^{n^2})$, with trace weights $2^{-n}$. Choose pairwise orthogonal projections $q_n$ in $\R$ of trace $2^{-n}$, $n\geq1$. For each $n$, the corner $q_n\R q_n$, with its normalized trace, contains a copy of $\bbM_n$ whose relative commutant contains a diffuse abelian von Neumann algebra with separable predual. Therefore $\bbM_n\tensor L_\infty((0,1)\times\bbT^{n^2})$ embeds into $q_n\R q_n$ by a trace-preserving $*$-homomorphism. Taking the direct sum over $n$ gives a trace-preserving embedding of $\N$ into $\R$. Thus, the space $E$ in Theorem~\ref{thm:main} can be chosen to satisfy
\[
 E\emb L_p(\R),\qquad E\nemb L_p(0,1).
\]
The application of \cite[Theorem~6]{JSZ} in the proof of Theorem~\ref{thm:main} gives an embedding of $C_E$ into the $L_p$-space of a finite von Neumann algebra, not necessarily into $L_p(\R)$.
\end{remark}

\begin{remark}
The random variables $Z_k$ defined in \eqref{eq:Z} are independent in the usual probabilistic sense as $X$-valued random variables. Under the embedding of $L_p(0,1;L_p(\A,\tau))$ into
\[
 L_p(L_\infty(0,1)\tensor\A,dm\otimes\tau),
\]
the images of $Z_k$ are not asserted to be independent relative to the scalar trace $dm\otimes\tau$. Consequently, the independence of the $X$-valued functions cannot be used to infer that the norm in \eqref{eq:E} admits a scalar-valued probabilistic representation of the kind studied in \cite{DC,RS}.
\end{remark}

\section{A continuous counterpart}

We record a consequence of Theorem~\ref{thm:main} and \cite[Corollary~2]{HJSZ} for symmetric spaces on the half-line. For a measurable scalar function $f$ on $(0,\infty)$, let $\mu(f)$ denote its decreasing rearrangement. Recall that the norm of $L_p+L_\infty$ on $(0,\infty)$ can be taken to be $(\int_0^1\mu(f)(s)^p\,ds)^{1/p}$.

\begin{corollary}\label{cor:continuous}
For every $1\leq p<2$, there exists a separable symmetric Banach function space $G(0,\infty)$ such that
\begin{enumerate}
\item[(i)] $G(0,\infty)\nemb L_p(0,1)$;
\item[(ii)] for every semifinite $\sigma$-finite von Neumann algebra $(\mathcal B,\nu)$, there exists a noncommutative probability space $(\widehat\N,\widehat\tau)$ such that $G(\mathcal B,\nu)\emb L_p(\widehat\N,\widehat\tau)$.
\end{enumerate}
Here $G(\mathcal B,\nu)$ is the symmetric operator space associated with $G(0,\infty)$ and $(\mathcal B,\nu)$.
\end{corollary}
\begin{proof}
Let $E$ be the symmetric sequence space constructed in the proof of Theorem~\ref{thm:main}. Following \cite[Section~1.1]{HJSZ}, define $G(0,\infty)$ to be the space of all $f\in(L_p+L_\infty)(0,\infty)$ such that the sequence
\[
 \left(\int_{k-1}^k\mu(f)(s)\,ds\right)_{k\geq1}
\]
belongs to $E$, equipped with the norm
\[
 \norm{f}_G=\left(\int_0^1\mu(f)(s)^p\,ds\right)^{1/p}
           +\left\|\left(\int_{k-1}^k\mu(f)(s)\,ds\right)_{k\geq1}\right\|_E.
\]
For $a=(a_k)_{k\geq1}\in E$, set $f_a=\sum_{k=1}^\infty a_k\ind_{(k-1,k)}$. Since $E\subset c_0$, the decreasing rearrangement of $f_a$ is the step function whose successive values are the coordinates of $\mu(a)$. Thus
\[
 \norm{f_a}_G=\norm{a}_\infty+\norm{a}_E,
 \qquad a\in E.
\]
In particular, $a\mapsto f_a$ is an isomorphic embedding of $E$ into $G(0,\infty)$, since $\norm{a}_\infty\leq\norm{a}_E$. Assertion (i) follows from $E\nemb L_p(0,1)$.

As shown in the proof of Theorem~\ref{thm:main}, $E\neq\ell_p$ and $E\emb L_p(\N,\rho)$ for a noncommutative probability space $(\N,\rho)$. Applying \cite[Corollary~2, (a)$\Rightarrow$(b)]{HJSZ}, we obtain an isomorphic embedding of $G(0,\infty)$ into the $L_p$-space of a finite von Neumann algebra. The same reference, namely \cite[Remark~10]{HJSZ}, gives the separability of $G(0,\infty)$. Assertion (ii) follows from the implication (b)$\Rightarrow$(c) in \cite[Corollary~2]{HJSZ}, including the extension to arbitrary semifinite $\sigma$-finite von Neumann algebras stated after that corollary.
\end{proof}

\end{document}